\documentclass{article}
\usepackage{amssymb}
\newtheorem{theorem}{Theorem}

\newtheorem{definition}[theorem]{Definition}
\newtheorem{example}[theorem]{Example}
\newtheorem{lemma}[theorem]{Lemma}

\newenvironment{proof}[1][Proof]{\textbf{#1.} }{\ \rule{0.5em}{0.5em}}

\title{Hausdorff Continuous Viscosity Solutions of Hamilton-Jacobi Equations}
\author{ Roumen Anguelov and Froduald Minani\\
Department of Mathematics and Applied Mathematics\\
University of Pretoria, Pretoria 0002\\
roumen.anguelov@up.ac.za\\
minanif2006@yahoo.fr}
\date{}
\begin{document}
\maketitle

\begin{abstract}
A new concept of viscosity solutions, namely, the Hausdorff
continuous viscosity solution for the Hamilton-Jacobi equation is
defined and investigated. It is shown that the main ideas within
the classical theory of continuous viscosity solutions can be
extended to the wider space of Hausdorff continuous functions
while also generalizing some of the existing concepts of
discontinuous solutions.
\end{abstract}



\section{Introduction}
The theory of viscosity solutions was developed for certain types
of first and second order PDEs. It has been particularly useful in
describing the solutions of PDEs associated with deterministic and
stochastic optimal control problems \cite{Fleming}, \cite{Bardi}.
In its classical formulation, see \cite{UsersGuide}, the theory
deals with solutions which are continuous functions. The concept
of continuous viscosity solutions was further generalized in
various ways, e.g. see [6, Chapter V], \cite{Barron},
\cite{Barles}, to include discontinuous solutions with the
definition of Ishii given in \cite{Ishii} playing a pivotal role.
In this paper we propose a new approach to the treatment of
discontinuous solutions, namely, by involving Hausdorff continuous
(H-continuous) interval valued functions. In the sequel we will
justify the advantages of the proposed approach by demonstrating
that
\begin{itemize}
\item the main ideas within the classical theory of continuous
viscosity solutions can be extended almost unchanged to the wider
space of H-continuous functions%
\item the existing theory of discontinuous solutions is a
particular case of that developed in this paper in terms of H-continuous functions%
\item the H-continuous viscosity solutions have a more clear
interpretation than the existing concepts of discontinuous
solutions, e.g. envelope viscosity solutions \cite[Chapter
V]{Bardi}.
\end{itemize}

In order to simplify the exposition we will only consider first
order Hamilton-Jacobi equations of the form
\begin{equation}\label{eq}
\Phi(x,u(x),Du(x))=0,\ x\in\Omega,
\end{equation}
where $\Omega$ is an open subset of $\mathbb{R}^n$,
$u:\Omega\rightarrow\mathbb{R}$ is the unknown function, $Du$ is
the gradient of $u$ and the given function
$\Phi:\Omega\times\mathbb{R}\times\mathbb{R}^n\rightarrow\mathbb{R}$
is jointly continuous in all its arguments.

The theory of viscosity solutions rests on two fundamental
concepts, namely, of subsolution and of supersolution. These
concepts are defined in various equivalent ways in the literature.
The definition given below is formulated in terms of local maxima
and minima. We will use the following notations
\begin{eqnarray*}
USC(\Omega)&=&\{u:\Omega\rightarrow\mathbb{R}:u\mbox{ is upper
semi-continuous on }\Omega\}\\
LSC(\Omega)&=&\{u:\Omega\rightarrow\mathbb{R}:u\mbox{ is lower
semi-continuous on }\Omega\}
\end{eqnarray*}
\begin{definition}\label{defsubsupsol}
A function $u\in USC(\Omega)$ is called a viscosity subsolution
of the equation (\ref{eq}) if for any $\varphi\in C^1(\Omega)$ we
have
\[
\Phi(x_0,u(x_0),D\varphi(x_0))\leq 0
\]
at any local maximum point $x_0$ of $u-\varphi$. Similarly, $u\in
LSC(\Omega)$ is called a viscosity supersolution of the equation
(\ref{eq}) if for any $\varphi\in C^1(\Omega)$ we have
\[
\Phi(x_0,u(x_0),D\varphi(x_0))\geq 0
\]
at any local minimum point $x_0$ of $u-\varphi$.
\end{definition}
Without loss of generality we may assume in the above definition
that $u(x_0)=\varphi(x_0)$ exposing in this way a very clear
geometrical meaning of this definition: the gradient of the
solution $u$ of equation (\ref{eq}) is replaced by the gradient of
any smooth function touching the graph of $u$ from above, in the
case of subsolution, and touching the graph of $u$ from below, in
the case of supersolution. This also establishes the significance
of the requirement that a subsolution and a supersolution should
respectively be upper semi-continuous and lower semi-continuous
functions. More precisely, the upper semi-continuity of a
subsolution $u$ ensures that any local supremum of $u-\varphi$ is
effectively reached at a certain point $x_0$, that is, it is a
local maximum, with the geometrical meaning that the graph of $u$
can be touched from above at $x=x_0$  by a vertical translate of
the graph of $\varphi$. In a similar way, the lower
semi-continuity of a supersolution $u$ ensures that any local
infimum of $u-\varphi$ is effectively reached at a certain point
$x_0$ which means that the graph of $u$ can be touched from below
at $x=x_0$ by a vertical translate of the graph of $\varphi$.

Naturally, a solution should be required somehow to incorporate
the properties of both a subsolution and a supersolution. In the
classical viscosity solutions theory, see \cite{UsersGuide}, a
viscosity solution is a function $u$ which is both a subsolution
and a supersolution. Since $USC(\Omega)\bigcap
LSC(\Omega)=C(\Omega)$, this clearly implies that the viscosity
solutions defined in this way are all continuous functions.

The concept of viscosity solution for functions which are not
necessarily continuous is introduced by using the upper and lower
semi-continuous envelopes, see \cite{Ishii}. Let us recall that
the upper semi-continuous envelope of a function $u$ which we
denote by $S(u)$ is the least upper semi-continuous function which
is not smaller than $u$. In a similar way, the lower
semi-continuous envelope $I(u)$ of a function $u$ is the largest
lower semi-continuous function not greater than $u$. For a locally
bounded function $u$ we have the following representations of
$S(u)$ and $I(u)$:
\begin{eqnarray}
S(u)(x)&=&\inf\{f(x):f\!\in\! USC(\Omega), u\leq f\}
=\inf_{\delta>0}\sup\{u(y):y\!\in\! B_\delta(x)\},\label{usce}\\
I(u)(x)&=&\sup\{f(x):f\!\in\! LSC(\Omega), u\geq f\}=
\sup_{\delta>0}\inf\{u(y):y\!\in\! B_\delta(x)\},\label{lsce}
\end{eqnarray}
where $B_\delta(x)$ denotes the open $\delta$-neighborhood of $x$
in $\Omega$. Using the fact that for any function
$u:\Omega\rightarrow\mathbb{R}$ the functions $S(u)$ and $I(u)$
are always, respectively, upper semi-continuous and lower
semi-continuous functions, a viscosity solution can be defined as
follows, \cite{Ishii}.

\begin{definition}\label{defvsol}
A function $u:\Omega\rightarrow\mathbb{R}$ is called a viscosity
solution of (\ref{eq}) if $S(u)$ is a viscosity subsolution of
(\ref{eq}) and $I(u)$ is a viscosity supersolution of (\ref{eq}).
\end{definition}

The first important point to note about the advantages of the
method in this paper is as follows. Interval valued functions
appear naturally in the context of noncontinuous viscosity
solutions. Namely, they appear as \emph{graph completions}.
Indeed, the above definition places requirements not on the
function $u$ itself but on its lower and upper semi-continuous
envelopes or, in other words, on the interval valued function
\begin{equation}\label{gc1}
F(u)(x)=[I(u)(x),S(u)(x)],\ x\in\Omega,
\end{equation}
which is called the \emph{graph completion} of $u$, see
\cite{Sendov}. Clearly, Definition \ref{defvsol} treats functions
which have the same upper and lower semi-continuous envelopes,
that is, have the same graph completion, as identical functions.
On the other hand, since different functions can have the same
graph completion, a function can not in general be identified from
its graph completion, that is, functions with the same graph
completion are indistinguishable. Therefore, no generality will be
lost if only interval valued functions representing graph
completions are considered.

Let $\mathbb{A}(\Omega )$ be the set of all functions defined on
an open set $\Omega \subset \mathbb{R}^{n}$ with values which are
closed finite real intervals, that is,
\[
\mathbb{A}(\Omega )=\{f:\Omega \rightarrow \mathbb{IR}\},
\]%
where $\mathbb{IR}=\{[\underline{a},\overline{a}]:\underline{a%
},\overline{a}\in \mathbb{R},\;%
\underline{a}\leq \overline{a}\}$. Identifying $a\in \mathbb{R}$
with the point interval $[a,a]\in \mathbb{IR}$, we consider
$\mathbb{R}$ as a subset of $\mathbb{IR}$. Thus $\mathbb{A}(\Omega
)$ contains the set $\mathcal(\Omega)=\{f:\Omega \rightarrow
\mathbb{R}\}$ of all real functions defined on $\Omega$.

Let $u\in \mathbb{A}(\Omega)$. For every $x\in \Omega$ the value
of $u$ is an interval $[\underline{u}(x),\overline{u}(x)]\in
\mathbb{IR}$. Hence, the function $u$ can be written in the form
$u=[\underline{u},\overline{u}]$ where $\underline{u},\overline{u
}\in \mathcal{A}(\Omega)$ and $\underline{u}(x)\leq
\overline{u}(x),\;x\in\Omega$. The function
\[
w(f)(x)=\overline{u}(x)-\underline{u}(x),\ x\in\Omega,
\]
is called width of $u$. Clearly, $u\in\mathcal{A}(\Omega)$ if and
only if $w(f)=0$. The definitions of the upper semi-continuous
envelope, the lower semi-continuous envelope and the graph
completion operator $F$ given in (\ref{usce}), (\ref{lsce}) and
(\ref{gc1}) for $u\in\mathcal{A}(\Omega)$ can be extended to
functions $u=[\underline{u},\overline{u}]\in\mathbb{A}(\Omega )$
as follows:
\begin{eqnarray*}
S(u)&=&\inf_{\delta>0}\sup\{z\in u(y):y\in B_\delta(x)\}\ =\
S(\overline{u})\\
I(u)&=&\sup_{\delta>0}\inf\{z\in u(y):y\in B_\delta(x)\}\ =\
I(\underline{u})\\
F(u)&=&[I(u),S(u)]\ =\ [I(\underline{u}),S(\overline{u})].
\end{eqnarray*}
We recall here the concept of S-continuity associated with the
graph completion operator, \cite{Sendov}.
\begin{definition}
A function $u=[\underline{u},\overline{u}]\in\mathbb{A}(\Omega)$
is called S-continuous if $F(u)=u$, or, equivalently,
$I(\underline{u})=\underline{u},\ S(\overline{u})=\overline{u}$.
\end{definition}
Using the properties of the lower and upper semi-continuous
envelopes one can easily see that the graph completions of locally
bounded real functions on $\Omega$ comprise the set
$\mathbb{F}(\Omega)$ of all S-continuous functions on $\Omega$.
Following the above discussion we define the concept of viscosity
solution for the interval valued functions in
$\mathbb{F}(\Omega)$.

\begin{definition}\label{defintvsol}
A function $u=[\underline{u},\overline{u}]\in\mathbb{F}(\Omega)$
is called a viscosity solution of (\ref{eq}) if $\underline{u}$ is
a supersolution of (\ref{eq}) and $\overline{u}$ is a subsolution
of (\ref{eq}).
\end{definition}

A second advantage of the method in this paper is as follows. A
function $u\in\mathcal{A}(\Omega)$ is a viscosity solution of
(\ref{eq}) in the sense of Definition \ref{defvsol} if and only if
the interval valued function $F(u)$ is a viscosity solution of
(\ref{eq}) in the sense of Definition \ref{defintvsol}. In this
way the level of the regularity of a solution $u$ is manifested
through the width of the interval valued function $F(u)$. It is
well known that without any additional restrictions the concept of
viscosity solution given in Definition \ref{defvsol} and by
implication the concept given in Definition \ref{defintvsol} is
rather weak, \cite{Bardi}. This is demonstrated by the following
example, which is also partially discussed in \cite{Bardi}.
\begin{example}\label{Ex1}
Consider the equation
\begin{equation}\label{Ex1eq}
u'(x)=1\ , \ \ x\in(0,1).
\end{equation}
The functions
\[
v(x)=\left\{\begin{tabular}{cll}$x\!+\!1$&if&$x\!\in\!(0,1)\!\cap\!\mathbb{Q}$\\
$x$&if&$x\!\in\!(0,1)\!\setminus\!\mathbb{Q}$\end{tabular}\right.
\ \ \ \
w(x)=\left\{\begin{tabular}{cll}$x$&if&$x\!\in\!(0,1)\!\cap\!\mathbb{Q}$\\
$x\!+\!1$&if&$x\!\in\!(0,1)\!\setminus\!\mathbb{Q}$\end{tabular}\right.
\]
are both viscosity solutions of equation (\ref{Ex1eq}) in terms of
Definition \ref{defvsol}. The interval valued function $
z=F(v)=F(w)$ given by
\[
z(x)=[x,x+1],\ x\in (0,1)
\]
is a solution in terms of Definition \ref{defintvsol}.
\end{example}

With the interval approach adopted here it becomes apparent that
the distance between $I(u)$ and $S(u)$ is an essential measure of
the regularity of any solution $u$, irrespective of whether it is
given as a point valued function or as an interval valued
function. If no restriction is placed on the distance between
$I(u)$ and $S(u)$ we will have some quite meaningless solutions
like the solutions in Example \ref{Ex1}. On the other hand, a
strong restriction like $I(u)=S(u)$ gives only solutions which are
continuous. In this paper we consider solutions for which the
Hausdorff distance, as defined in \cite{Sendov}, between the
functions $I(u)$ and $S(u)$ is zero, a condition defined through
the concept of Hausdorff continuity.

\section{The space of Hausdorff continuous functions}

The concept of Hausdorff continuous interval valued functions was
originally developed within the theory of Hausdorff
approximations, \cite{Sendov}. It generalizes the concept of
continuity of real function using a minimality condition with
respect to inclusion of graphs.

\begin{definition}
\label{dhcon}A function $f\in\mathbb{A}(\Omega)$ is called
Hausdorff continuous, or H-continuous, if for every
$g\in\mathbb{A}(\Omega)$ which satisfies the inclusion
$g(x)\subseteq f(x),$ $x\in \Omega,$ we have $F(g)(x)=f(x),$ $x\in
\Omega$.
\end{definition}

The following theorem gives useful necessary and sufficient
conditions for an interval valued function to be H-continuous,
\cite{Sendov}, \cite{Anguelov}.
\begin{theorem}\label{thcondi}
Let $f=[\underline{f},\overline{f}]\in\mathbb{A}(\Omega)$. The
following conditions are equivalent
\begin{itemize}
\item[a)] the function $f$ is H-continuous
\item[b)]$F(\underline{f})=F(\overline{f})=f$%
\item[c)] $S(\underline{f})=\overline{f}$,
$I(\overline{f})=\underline{f}$ and $f$ is S-continuous
\end{itemize}
\end{theorem}

As mentioned in the Introduction the concept of Hausdorff
continuity is closely connected with the Hausdorff distance
between functions as introduced by Sendov in \cite{Sendov}. The
Hausdorff distance $\rho(f,g)$ between two functions
$f,g\in\mathbb{A}(\Omega)$ is defined as the Hausdorff distance
between the graphs of the functions $F(f)$ and $F(g)$ considered
as subsets of $\mathbb{R}^{n+1}$. More precisely  we have
\begin{eqnarray*}
\rho(f,g)&=&\max\{\sup_{x_1\in\Omega}\sup_{y_1\in
F(f)(x_1)}\inf_{x_2\in\Omega}\inf_{y_2\in
F(g)(x_2)}||(x_1-x_2,y_1-y_2)||,\\
&&\sup_{x_2\in\Omega}\sup_{y_2\in
F(g)(x_2)}\inf_{x_1\in\Omega}\inf_{y_1\in
F(f)(x_1)}||(x_1-x_2,y_1-y_2)||\}\ .
\end{eqnarray*}
where $||\cdot||$ is a given norm in $\mathbb{R}^{n+1}$. Condition
b) in the Theorem \ref{thcondi} implies that for any H-continuous
function $f=[\underline{f},\overline{f}]$ the Hausdorff distance
between the functions $\underline{f}$ and $\overline{f}$ is zero.
More precisely we have
\[f=[\underline{f},\overline{f}]\ \mbox{ is H-continuous}\ \Longleftrightarrow \
\left\{\begin{tabular}{l}$f$ is
S-continuous\\\\$\rho(\underline{f},\overline{f})=0$.\end{tabular}\right.
\]

Although every H-continuous function $f$ is, in general, interval
valued, the subset of the domain $\Omega$ where $f$ assumes proper
interval values is a set of first Baire category. This result is
stated in the following theorem where it is also shown that for
H-continuous functions interval values are used in an 'economical'
way, namely only at points of discontinuity, \cite{Anguelov}.

\begin{theorem}
\label{tcont}Let $f=[\underline{f},\overline{f}]$ be an
H-continuous function on $\Omega$.%
\begin{itemize}\item[a)] If $\underline{f}$ or $\overline{f}$ is continuous at a point $a\in
\Omega$ then $\underline{f}(a)=\overline{f}(a)$. %
\item[b)] If $\underline{f}(a)=\overline{f}(a)$ for some $a\in \Omega$ then both $%
\underline{f}$ and $\overline{f}$ are continuous at $a$. %
\item[c)]The set
\[
W_f=\{x\in\Omega : w(f(x))>0\}
\]
is a set of first Baire category.
\end{itemize}
\end{theorem}

Further properties of the H-continuous functions are discussed in
\cite{Sendov}, \cite{Anguelov-Markov}, \cite{Anguelov}, where it
is shown, among others, that they retain some of the essential
characteristics of the usual continuous functions. For example, an
H-continuous function is completely determined by its values on
any dense subset of the domain as stated in the following theorem
\cite{Anguelov}:
\begin{theorem} \label{tindent}
Let $f,g$ be H-continuous on $\Omega$ and let $D$ be a dense subset of $\Omega$. Then%
\begin{itemize}
\item[a)] $\ f(x)\ \leq\  g(x),\ x\in D\ \Longrightarrow\ \
f(x)\leq g(x),\ x\in \Omega,$ \item[b)] $\ f(x)\  =\ g(x),\ x\in
D\ \Longrightarrow\ f(x)=g(x),\ x\in \Omega.$
\end{itemize}
\end{theorem}

One of the most surprising and useful properties of the set
$\mathbb{H}(\Omega)$ of all H-continuous functions is its Dedekind
order completeness. What makes this property so significant is the
fact that with very few exceptions the usual spaces in Real
Analysis or Functional Analysis are not Dedekind order complete.
The order considered in $\mathbb{H}(\Omega)$ is the one which is
introduced point-wise, \cite{Anguelov}, \cite{Markov}, as follows:
For $f=[\underline{f},\overline{f}]\in\mathbb{H}(\Omega)$ and
$g=[\underline{g},\overline{g}]\in\mathbb{H}(\Omega)$ we have
\begin{equation}\label{forder}
 f\leq g\Longleftrightarrow \underline{f}(x)\leq
\underline{g}(x),\ \overline{f}(x)\leq \overline{g}(x),\
x\in\Omega.
\end{equation}
\begin{theorem}\label{tocomp}
The set $\mathbb{H}(\Omega)$ of all H-continuous interval valued
functions is Dedekind order complete with respect to the order
defined through (\ref{forder}), that is,%
\begin{itemize}\item[(i)] for every subset
$\mathcal{F}$ of $\, \mathbb{H}(\Omega)$ which is bounded from
above there exist $u\in\mathbb{H}(\Omega)$ such that
$u=\sup\mathcal{F}$ %
\item[(ii)]for every subset $\mathcal{F}$ of $\,
\mathbb{H}(\Omega)$ which is bounded from below there exist
$v\in\mathbb{H}(\Omega)$ such that $v=\inf \mathcal{F}$.
\end{itemize}
\end{theorem}

We should note that the supremum and infimum in the above theorem
are not defined in a point-wise way and that the point-wise
supremum and infimum and not necessarily H-continuous functions.
The following representation of the supremum in the poset
$\mathbb{H}(\Omega)$ through the point-wise supremum is useful,
\cite{Anguelov2}.

\begin{theorem}\label{tpointsup}
Let the set $\mathcal{F}\subseteq\mathbb{H}(\Omega)$ be bounded
from above and let the function $\psi\in\mathcal{A}(\Omega)$ be defined by%
\[
\psi(x)=\sup\{\overline{f}(x):f=[\underline{f},\overline{f}]\in\mathcal{F}\}\
,{\rm \ }x\in \Omega.
\]
Then
\[
\sup\mathcal{F}=F(S(\psi)).
\]
\end{theorem}
A similar representation holds for the infimum in the set
$\mathbb{H}(\Omega)$.

\section{The envelope viscosity solutions and  Hausdorff continuous
viscosity solutions}

Recognizing that the concept of viscosity solution given by
Definition \ref{defvsol} is rather weak the authors of
\cite{Bardi} introduce the concept of envelope viscosity solution.
The concept is defined in \cite{Bardi} for the equation (\ref{eq})
with Dirichlet boundary conditions. In order to keep the
exposition as general as possible we will give the definition
without explicitly involving the boundary condition.
\begin{definition}\label{desol}
A function $u\in\mathcal{A}(\Omega)$ is called an envelope
viscosity solution of (\ref{eq}) if there exist a nonempty set
$\mathcal{Z}_1(u)$ of subsolutions of (\ref{eq}) and a nonempty
set $\mathcal{Z}_2(u)$ of supersolutions of (\ref{eq}) such that
\[
u(x)=\sup_{f\in\mathcal{Z}_1(u)}f(x)=\inf_{f\in\mathcal{Z}_2(u)}f(x),\
x\in\Omega.
\]
\end{definition}
It is shown in \cite{Bardi} that every envelope viscosity solution
is a viscosity solution in terms of Definition \ref{defvsol}.
Considering the concept from geometrical point of view, on can
expect that by 'squeezing' the envelope viscosity solution $u$
between a set of subsolutions and a set of supersolutions the gap
between $I(u)$ and $S(u)$ would be small. However, in general this
is not the case. The following example shows that the concept of
envelope viscosity solution does not address the problem of the
distance between $I(u)$ and $S(u)$. Hence one can have envelope
viscosity solutions of little practical meaning similar to the
viscosity solution in Example \ref{Ex1}.

\begin{example}\label{Ex2} Consider the following equation on
$\Omega=(0,1)$
\begin{equation}\label{Ex2eq}
-u(x)(u'(x))^2=0,\ x\in \Omega.
\end{equation}
For every $\alpha\in\Omega$ we define the functions
\begin{eqnarray*}
\phi_\alpha(x)&=&\left\{\begin{tabular}{l}$1$ if $x=\alpha$\\ $0$
if $x\in\Omega\setminus\{\alpha\}$\end{tabular}\right.\\
\psi_\alpha(x)&=&\left\{\begin{tabular}{l}$0$ if $x=\alpha$\\ $1$
$x\in\Omega\setminus\{\alpha\}$.\end{tabular}\right.
\end{eqnarray*}
We have
\[
\phi_\alpha\in USC(\Omega), \ \psi_\alpha\in LSC(\Omega),\ \alpha\in\Omega.
\]
Furthermore, for every $\alpha\in (0,1)$ the function
$\phi_\alpha$ is a subsolution of (\ref{Ex2eq}) while
$\psi_\alpha$ is a supersolution of (\ref{Ex2eq}). Indeed, both
functions satisfy the equation for all $x\in
\Omega\setminus\{\alpha\}$ and at $x=\alpha$ we have
\begin{eqnarray*}
-\phi_\alpha(\alpha)p^2&=&-p^2 \leq 0\ \mbox{ for all }\ p\in
D^+\phi_\alpha(\alpha)=(-\infty,\infty)\\
-\psi_\alpha(\alpha)p^2&=&0\geq 0\ \mbox{ for all }\ p\in
D^-\psi_\alpha(\alpha)=(-\infty,\infty).
\end{eqnarray*}
We will show that the function
\[
u(x)=\left\{\begin{tabular}{l}$1$ if
$x\in\Omega\setminus\mathbb{Q}$ \\ $0$ if
$x\in\mathbb{Q}\bigcap\Omega$\end{tabular}\right.
\]
is an envelope viscosity solution of (\ref{Ex2eq}). Define
\begin{eqnarray*}
\mathcal{Z}_1(u)&=&\{\phi_\alpha:\alpha\in\Omega\setminus\mathbb{Q}\}\\
\mathcal{Z}_2(u)&=&\{\psi_\alpha:\alpha\in\mathbb{Q}\bigcap\Omega\}
\end{eqnarray*}
Then $u$ satisfies
\[
u(x)=\sup_{w\in\mathcal{Z}_1(u)} w(x)=\inf_{w\in\mathcal{Z}_2(u)}
w(x)
\]
which implies that it is an envelope viscosity solution. Clearly
neither $u$ nor $F(u)$ is a Hausdorff continuous function. In fact
we have $F(u)(x)=[0,1]$, $x\in\Omega$.
\end{example}

The next interesting question is whether every H-continuous
solution is an envelope viscosity solution. Since the concept of
envelope viscosity solutions requires the existence of sets of
subsolutions and supersolutions respectively below and above an
envelope viscosity solution then an H-continuous viscosity
solution is not in general an envelope viscosity solution, e.g.
when the H-continuous viscosity solutions does not have any other
subsolutions and supersolutions around it. However in the
essential case when the H-continuous viscosity solution is a
supremum of subsolutions or infimum of supersolutions it can be
linked to an envelope viscosity solution as stated in the next
theorem.

\begin{theorem}
Let $u=[\underline{u},\overline{u}]$ be an H-continuous viscosity
solution of (\ref{eq}) and let
\begin{eqnarray*}
Z_1&=&\{w\in USC(\Omega):w\mbox{-subsolution},\
w\leq\underline{u}\}\\
Z_2&=&\{w\in LSC(\Omega):w\mbox{-supersolution},\
w\geq\overline{u}\}.
\end{eqnarray*}
\begin{itemize}
\item[a)] If $Z_1\neq\emptyset$ and
$\underline{u}(x)=\sup\limits_{w\in Z_1}w(x)$ then $\underline{u}$
is an envelope viscosity solution.

\item[b)] If $Z_2\neq\emptyset$ and
$\overline{u}(x)=\inf\limits_{w\in Z_2}w(x)$ then $\overline{u}$
is an envelope viscosity solution.
\end{itemize}
\end{theorem}
\begin{proof}
a) We choose the sets $\mathcal{Z}_1(\underline{u})$ and
$\mathcal{Z}_2(\underline{u})$ required in Definition \ref{desol}
as follows
\[
\mathcal{Z}_1(\underline{u})=Z_1\ ,\ \
\mathcal{Z}_2(\underline{u})=\{\underline{u}\}.
\]
Then we have
\[
\underline{u}(x)=\sup_{w\in\mathcal{Z}_1(\underline{u})}w(x)=\inf_{w\in\mathcal{Z}_2(\underline{u})}w(x)
\]
which implies that $\underline{u}$ is an envelope viscosity
solution.

The proof of b) is done in a similar way.
\end{proof}

Let us note that if the conditions in both a) and b) in the above
theorem are satisfied then both $\underline{u}$ and $\overline{u}$
are envelope viscosity solutions and in this case in makes even
more sense to consider instead the H-continuous function $u$.

\section{Existence of Hausdorff continuous
viscosity solutions}

One of the primary virtues of the theory of viscosity solutions is
that it provides very general existence and uniqueness theorems,
\cite{UsersGuide}. In this section we will formulate and prove an
existence theorems for H-continuous viscosity solutions in a
similar form to the respective theorems for continuous solutions,
\cite{UsersGuide} (Theorem 4.1), and for general discontinuous
solutions, \cite{Ishii} (Theorem 3.1), \cite{Bardi} (Theorem
V.2.14).

\begin{theorem}\label{texist}
Assume that there exists Hausdorff continuous functions
$u_1=[\underline{u}_1,\overline{u}_1]$ and
$u_2=[\underline{u}_2,\overline{u}_2]$ such that $\overline{u}_1$
is a subsolution of (\ref{eq}), $\underline{u}_2$ is a
supersolution of (\ref{eq}) and $u_1\leq u_2$. Then there exists a
Hausdorff continuous solution $u$ of (\ref{eq}) satisfying the
inequalities
\[
u_1\leq u\leq u_2.
\]
\end{theorem}
The proof of the above theorem, similar to the other existence
theorems in the theory of viscosity solutions, uses Perron's
method and the solutions will be constructed as a supremum of a
set of subsolutions, this time the supremum being taken in the
poset $\mathbb{H}(\Omega)$ and not point-wise. We should note that
due to the fact that the the poset $\mathbb{H}(\Omega)$ is
Dedekind order complete it is an appropriate medium for such an
application of Perron's method. In the proof we will also use the
so called 'Bump Lemma', which can be formulated for Hausdorff
continuous functions as follows.

\begin{lemma}\label{lbump}
Let $u=[\underline{u},\overline{u}]\in\mathbb{H}(\Omega)$ be such
that $\overline{u}$ is a subsolution of (\ref{eq}) and
$\underline{u}$ fails to be a supersolution of (\ref{eq}) at some
point $y\in\Omega$. Then, for any $\delta>0$ there exists
$\gamma>0$ such that, for all $r<\gamma$, there exists a function
$w=[\underline{w},\overline{w}]\in\mathbb{H}(\Omega)$ with the
following properties:
\begin{itemize}
\item[(i)]$\overline{w}$ is a subsolution of (\ref{eq}),%
\item[(ii)]$w\geq u$,%
\item[(iii)]$w\neq u$,%
\item[(iv)]$w(x)=u(x),\ x\in\Omega\setminus B_r(y)$,%
\item[(v)]$\underline{w}(x)\leq \max\{\underline{u}(x),
\underline{u}(y)+\delta\},\ x\in B_r(y)$.
\end{itemize}
\end{lemma}

The proof of Lemma \ref{lbump} is similar to the proof of the Bump
Lemma in \cite{Bardi} (Lemma V.2.12) for real function with some
obvious changes due to interval character of the functions $u$ and
$w$.

We will also use the following result which was proved in
\cite[Proposition V.2.11]{Bardi}.

\begin{theorem}\label{tV-2-11}
\begin{itemize}
\item[a)] Let $\mathcal{Z}_1\subseteq USC(\Omega)$ be a set of
subsolutions of (\ref{eq}). If the function
\[
u(x)=\sup_{w\in\mathcal{Z}_1}w(x),\ x\in\Omega
\]
is locally bounded then $S(u)$ is a subsolutions of (\ref{eq}).%
\item[b)] Let $\mathcal{Z}_2\subseteq LSC(\Omega)$ be a set of
supersolutions of (\ref{eq}). If the function
\[
v(x)=\inf_{w\in\mathcal{Z}_2}w(x),\ x\in\Omega
\]
is locally bounded then $I(v)$ is a supersolution of (\ref{eq}).
\end{itemize}
\end{theorem}

\noindent{\bf Proof of Theorem \ref{texist}.} Consider the set
\[
\mathcal{U}=\{w=[\underline{w},\overline{w}]\in\mathbb{H}(\Omega):w\leq
u_2,\ \overline{w}\mbox{ is a subsolution} \}\ .
\]
Clearly the set $\mathcal{U}$ is not empty since
$u_1\in\mathcal{U}$. Let $u=\sup\mathcal{U}$ where the supremum is
taken in the set $\mathbb{H}(\Omega)$, i.e.,
$u\in\mathbb{H}(\Omega)$. We will show that $u$ is the required
viscosity solution of (\ref{eq}). Obviously we have the
inequalities
\[
u_1\leq u\leq u_2.
\]
Furthermore, according to Theorem \ref{tpointsup}, $u$ is given by
\[
u=F(S(\psi))
\]
where
\[
\psi(x)=\sup\{\overline{w}(x):w=[\underline{w},\overline{w}]\in\mathcal{U}\},\
x\in\Omega.
\]
Using that $\overline{w}$ is a subsolution for all
$w=[\underline{w},\overline{w}]\in\mathcal{U}$, it follows from
Theorem \ref{tV-2-11} that $\overline{u}=S(\psi)$ is a
subsolution. It remains to show that $\underline{u}$ is a
supersolution. To this end let us fix $y\in\Omega$.

Consider first the case when
$\underline{u}(y)=\underline{u}_2(y)$. Let $\varphi\in
C^1(\Omega)$ be such that $\underline{u}-\varphi$ has a local
minimum at $y$ and $\underline{u}(y)=\varphi(y)$. Then, in a
neighborhood of $y$, we have
\[
(\underline{u}_2-\varphi)(x)\geq (\underline{u}-\varphi)(x)\geq
0=(\underline{u}_2-\varphi)(y).
\]
Therefore, the function $\underline{u}_2-\varphi$ also has a local
minimum at $y$. Using that $\underline{u}_2$ is a supersolution we
obtain
\[
\Phi(y,\underline{u}_2(y),D\varphi(y))\geq 0.
\]
Since $\underline{u}(y)=\underline{u}_2(y)$ the above inequality
shows that the function $\underline{u}$ satisfies at the point $y$
the conditions of supersolution as stated in Definition
\ref{defsubsupsol}.

Consider now the case when $\underline{u}(y)<\underline{u}_2(y)$.
Then there exists $\delta>0$ such that%
\begin{equation}\label{ineq1}
\underline{u}(y)+\delta\leq\underline{u}_2(y)-\delta.
\end{equation}
Assume that $\underline{u}$ fails to be a supersolution at the
point $y$. Then, according to Lemma \ref{lbump}, there exists a
function $w\in\mathbb{H}(\Omega)$ with the properties (i)-(v),
where, using also the lower semi-continuity of $\underline{u}_2$,
$r>0$ is chosen in such a way that
\begin{equation}\label{ineq2}
\underline{u}_2(y)-\delta\leq\underline{u}_2(x),\ x\in B_r(y).
\end{equation}
Using (\ref{ineq1}) and (\ref{ineq2}), we obtain
\[
\underline{u}(y)+\delta\leq\underline{u}_2(y)-\delta\leq
\underline{u}_2(x),\ x\in B_r(y).
\]
Hence, from property (v) of Lemma \ref{lbump} for $x\in B_r(y)$ we
have
\[
\underline{w}(x)\leq\max\{\underline{u}(x),\underline{u}(y)+\delta\}\leq\underline{u}_2(x)
\]
Due to property (iv) the above inequality can be extended to all
$x\in\Omega$ and we have $\underline{w}\leq\underline{u}_2$. Using
Theorem \ref{thcondi}b) this inequality can be transferred over to
the functions $w$ and $u_2$ as follows
\[
w=F(\underline{w})\leq F(\underline{u}_2)=u_2
\]
This implies that $w\in\mathcal{U}$. Then $u=\sup\mathcal{U}\geq
w$ which contradicts conditions (ii) and (iii) in Lemma
\ref{lbump}. The obtained contradiction shows that $\underline{u}$
is a supersolution. Therefore the H-continuous function $u$ is a
viscosity solution of (\ref{eq}) in terms of Definition
\ref{defintvsol}. \rule{0.5em}{0.5em}

\section{Conclusion}
The Hausdorff continuous functions, being a particular class of
interval valued functions, belong to what is usually called Interval
Analysis, see \cite{Moore}. Nevertheless, recent results have shown
that they can provide solutions to problems formulated in terms of
point valued functions. A long outstanding problem related to the
Dedekind order completion of spaces $C(X)$ of real valued continuous
functions on rather arbitrary topological spaces $X$ was solved
through Hausdorff continuous functions, \cite{Anguelov}. Following
this breakthrough a significant improvement of the regularity
properties of the solutions obtained through the order completion
method, see \cite{Rosinger}, was reported in \cite{AnguelovRosinger}
and \cite{AngRos1}. Namely, it was shown that these solutions can be
assimilated with the class of Hausdorff continuous functions on the
open domains $\Omega$.

In this paper the Hausdorff continuous functions are linked with
the concept of viscosity solutions. As shown in the Introduction
the definition of viscosity solution, see Definition
\ref{defvsol}, has an implicit interval character since it places
requirements only on the upper semi-continuous envelope $S(u)$ and
the lower semi-continuous envelope $I(u)$. For a Hausdorff
continuous viscosity solution $u$ the functions $I(u)$ and $S(u)$
are as close as they can be in the sense of the Hausdorff distance
$\rho$ defined in \cite{Sendov}, namely, we have
$\rho(I(u),S(u))=0$. Hence, the requirement that a viscosity
solution is Hausdorff continuous has a direct interpretation which
we find clearer than the requirements related to some other
concepts of discontinuous viscosity solutions. The first main
result in the paper is that the concept of envelope viscosity
solution, which is generally used to single out the physically
meaningful solutions, is a particular case of the concept of
Hausdorff continuous viscosity solution. The second main result,
an existence theorem for Hausdorff continuous solutions, shows
that the main ideas of the classical theory of viscosity solutions
can be extended to Hausdorff continuous solutions. Further
research will seek a suitable formulation of comparison principle
for Hausdorff continuous viscosity solutions and respective
uniqueness results.

\end{document}